\documentclass[a4paper,reqno]{amsart}

\usepackage{tgschola}

\usepackage{latexsym}
\usepackage[english]{babel}
\usepackage{fancyhdr}
\usepackage[mathscr]{eucal}
\usepackage{amsmath}
\usepackage{mathrsfs}
\usepackage{mathtools}
\usepackage{amsthm}
\usepackage{amsfonts}
\usepackage{amssymb}
\usepackage{amscd}
\usepackage{bbm}
\usepackage{graphicx}
\usepackage{graphics}
\usepackage{latexsym}
\usepackage{color}
\usepackage{pifont}
\usepackage{tikz}
\usepackage[normalem]{ulem}

\usepackage{geometry}

\theoremstyle{plain}

\theoremstyle{definition}

\newtheorem*{remark*}{Remark}

\numberwithin{equation}{section}

\begin{document}

\title{Singularity: a Seventh Memo?}

\author[M. Gallone]{Matteo Gallone}
\address[M.~Gallone]{SISSA (Scuola Internazionale Superiore di Studi Avanzati)\\ via Bonomea 265 \\ 34136 Trieste (Italy).}
\email{matteo.gallone@sissa.it}

\author[S. Lucente]{Sandra Lucente}
\address[S.~Lucente]{ Dipartimento Interateneo di Fisica
Università degli Studi di Bari Aldo Moro, via Amendola 173 - 70125 Bari, Italy} 
\email{sandra.lucente@uniba.it}

\date{\today}

\subjclass[2000]{}

%
%

\begin{abstract}{In this paper we explore the relationships between Calvino's memos and Mathematics. In the first part, we discuss how Lightness, Quickness, Exactitude, Visibility, Multiplicity are present in the mathematical language, reasoning and in the work of the mathematician. In addiction, we follow a similar path for the topics of Calvino's lecture of which we only know the title or some notes. In the final part, we explain why `Singularity' could be chosen as topic for a possible Calvino's seventh lecture.}
\end{abstract}
\maketitle

\tableofcontents

\section{Introduction}\label{sec:Introduction}
In 1985, Italo Calvino was to have held a series of lectures at Harvard entitled \emph{Six Memos for the Next Millennium}, but a fatal stroke left only five transparent folders in a hard folder on his desk. Each folder contained the text of a lecture with the following titles: \emph{Lightness}, \emph{Quickness}, \emph{Exactitude}, \emph{Visibility}, and \emph{Multiplicity}. In the Introduction of the Italian version \cite{C1}, we find that from Calvino's point of view the memos he discusses are ``certain values, qualities, or peculiarities of Literature that are very close to my heart.''  But are these values only literary? What would be the qualities of science to be developed in this era? In particular, can we seek specific characteristics of those who research in Mathematics nowadays?

In the Preface that the author's wife writes to accompany the posthumously published Italian text titled \emph{American Lessons}, Esther Judith Singer Calvino tells us about the period of preparation for the lectures: ``They became an obsession, and one day he announced to me that he had ideas and material for eight lectures. I know the title for what might have been an eighth lecture: \emph{Sul cominciare e sul finire} (On the beginning and the ending [of novels]).''  Notes are found for this lecture. It is also known that the sixth lecture would have been titled \emph{Consistency}. However, there is no reference to the seventh lecture.

As part of the \emph{Puglia Summer Trimester 2023}, the authors of this article held a lecture on the ``American Lessons''  and, recognizing in the sub-title of the trimester ``Singularities, Asymptotics, and Limiting Models''  three mathematical memos, they sought others in Calvino's text. 

This article is divided into three parts. Sections \ref{sec:Lightness}-\ref{sec:Multiplicity} explore the relationships between Mathematics and the five memos chosen by Calvino: Lightness, Quickness, Exactitude, Visibility, Multiplicity. Sections \ref{sec:Consistency}-\ref{sec:BeginEnd} try to follow a similar path for the lessons of which we only know the titles or some notes. In Section \ref{sec:Singularity}, we allow ourselves to play with the absence of the seventh lesson by choosing for it one of the memos from the conference: Singularity.

This attempt to apply to science the values of the American Lessons is not a new idea. The most successful mirror between the methods of Literature proposed by Calvino and that of Mathematics is the book by Gabriele Lolli \cite{L1} where a different vision of our paragraphs \ref{sec:Lightness}-\ref{sec:Multiplicity} can be found. Regarding paragraph \ref{sec:Singularity}, i.e., the creation of a seventh scientific lesson in Calvino's style, in 2010 Nicola Cabibbo held a lecture at the Accademia dei Lincei titled ``Science as the search for clarity''  and indeed he imagines a seventh Calvino-style lesson on ``Clarity''  (see \cite{C}). In this article, therefore, we do not want to perform any mathematical exegesis of the American Lessons, we only take up Calvino's invitation to be protagonist readers of the text with our characteristic of being mathematicians. This invitation in ``If on a winter's night a traveler'' \cite{C3} becomes imperative, the reader becomes Reader, here the mathematical reader becomes Author of a paper. Perhaps we will remain on the surface of Mathematics and Literature, but as Palomar says

\begin{quote}
Only after knowing the surface of things can one venture to seek what is beneath. But the surface of things is inexhaustible. \cite{C4}.
\end{quote}

\section{The Lightness of Mathematical Models}\label{sec:Lightness}

Calvino chooses \emph{Lightness} as the first of the ``Memos''  and we do the same in this first paragraph, showing that is central to Mathematics. This value is an essential characteristic that leads to the birth of all mathematical entities. In other words, we would like to convince the reader that the objects with which Mathematics works are \emph{light} in Calvino's sense. In the cases where Mathematics is the language of an applied science, such as Physics, there is an unavoidable lightening during the process of mathematization of reality through the construction of a model. As an example, let us look at the problem of the motion of celestial bodies, ingeniously solved by Newton with the law of universal gravitation. Today we know that the simplest way to describe the motion of the planets in the Solar System is in a reference frame where the Sun is fixed and the planets orbit around it. But let us pretend we don't know Newton's answer and put ourselves in the shoes of a supposed astronomer from the past aware of the fact that the Earth is a planet. If we want to describe the motion of the Earth around the Sun, we would probably have to start with the description of the Earth itself. We would then realize that the Earth is extremely difficult to describe: there are oceans and continents; on the continents, there are flat areas, mountains, rivers, and lakes. The flat areas can then be green, arid, or icy, and so on. In the attempt to construct the most accurate model of Earth, we would end up getting lost in the wonderful details of the planet we live in, losing sight of the goal we had set at the beginning of our investigation. In a somewhat colloquial jargon, we might comment by saying that this model, with all this wealth of details, is \emph{heavy} to handle. This is where the first of the quotations from Calvino's \emph{Memos} comes into play:
\begin{quote}
My working method has more often than not involved the subtraction of weight. I have tried to remove weight, sometimes from people, sometimes from heavenly bodies, sometimes from cities; above all I have tried to remove weight from the structure of stories and from language. \cite{C2}.
\end{quote}

The path suggested by reading Calvino is to remove weight, that is, to deprive the model of superfluous details. But what are the superfluous details? The question is legitimate and it is impossible to give a general answer. To give an idea of what guidelines to follow, Calvino quotes of the myth of Perseus and Medusa.

Briefly, Medusa is a Gorgon, that is, a woman with a terrifying appearance with snakes instead of hair who can petrify instantly anyone who crosses her gaze. The myth tells that a young hero named Perseus had to bring as a bridal gift the head of Medusa. To kill the Gorgon, he uses a very ingenious strategy and, instead of confronting her openly, prefers to use a highly polished shield to observe the reflected image of the Gorgon. After beheading Medusa, Perseus collects her head and carries it in a bag with him, making it a very powerful weapon to be used only ``with enemies who deserve to become the statue of themselves.'' \footnote{Here, we prefer to translate the Italian version. In the English version it reads: ``against those who deserve the punishment of being turned into statues''.} It is interesting to note that Calvino represents in the fight between Perseus and Medusa a clash between the \emph{lightness} (Perseus arrives flying on Pegasus, a winged horse) and the \emph{heavyness} (Medusa petrifies with her gaze). The clash is won by the Lightness embodied by Perseus, who carries with him a sign of heaviness, represented by the decapitated head of the Gorgon which gives him the ability to petrify, that is, to make heavy. To quote Calvino:

\begin{quote}
Perseus succeeds in mastering that horrendous face by keeping it hidden, just as in the first place he vanquished it by viewing it in a mirror. Perseus's strength always lies in a refusal to look directly, but not in a refusal of the reality in which he is fated to live; he carries the reality with him and accepts it as his particular burden. \cite{C2}.
\end{quote}

This description fits perfectly when looking at the work of lightening a model. In fact, every detail that is decided to be neglected to lighten the model, could actually be an essential one. Out of metaphor, what is said in this passage is that only by carrying with us in reasoning all the details that we have neglected we can perform the operation of lightening in a correct way. For this reason, Calvino quotes Ovid to emphasize with how much care these superfluous details should be treated:

\begin{quote}
So that the rough sand should not harm the snake-haired head he makes the ground soft with a bed of leaves, and on top of that he strews little branches of plants born under water, and on this he places Medusa's head, face down. \cite{C2}.
\end{quote}

Let's return to the operation of simplification in the case of analyzing the Earth's motion around the Sun. First, we eliminate all elevations, all seas, and all land masses, turning our planet into a sphere. Then, we reduce the size of this sphere until it becomes a point. Reducing the Earth to a dimensionless mathematical object, such as a point, is sufficient to fully understand its orbit around the Sun. Surprisingly, this simplification, at first glance, seems to be a ``loss of information''  because of the many neglected details. In fact, not only can the Earth be idealized as a point, but every planet in the Solar System can be idealized as a point. Therefore, in this context, in the simplification where it seemed we had lost information, we actually gain it. We gain that the laws we write for the motion of the Earth are the same that we would write for the motion of any other planet.

A similar discussion can be had for many objects in Mathematics. Abstract geometric entities capture the essence of a vast and diverse set of objects and situations from our daily lives. Circles, triangles, lines, as well as functions and more abstract objects like Banach spaces or differential manifolds arise from the ``simplification'' of reality in the sense of Calvino. After the simplification process, each such object becomes more versatile and general than the original ones. The relationships established between them form an intricate web of Lemmas, Propositions, and Theorems that enter into Calvino's ``thoughtful''  Lightness:

\begin{quote}
There is such a thing as a lightness of thoughtfulness, just as we all know that there is a lightness of frivolity. In fact, thoughtful lightness can make frivolity seem dull and heavy. \cite{C2}.
\end{quote}

\section{Geodetic and Zigzagging Quickness}\label{sec:Quickness}

It is surprising to think that Calvino chooses \emph{Quickness} as a value for Literature of the new millennium, and even more surprising is the idea of rapid Mathematics. What is too rapid today is the time dedicated to reading texts and even scientific articles. In the last ten years, 2.5 million scientific articles have been produced; to read them at a rate of one line per minute without interruptions would take the same amount of time that \emph{Homo Sapiens} has lived on this planet. We do not read; we scroll. And here Calvino helps us understand that sometimes this is a value. 

\begin{quote}
``Discoursing,''  or ``discourse,''  for Galileo means reasoning, and very often deductive reasoning. ``Discoursing is like coursing'': this statement could be Galileo's declaration of faith—style as acmethod of thought and as literary taste. For him, good thinking means Quickness, agility increasoning, economy in argument, but also the use of imaginative examples. \cite{C2}.
\end{quote}

For Calvino, Galileo is the greatest writer of the Italian Literature of every century. He wrote this clearly in the \emph{Corriere della Sera} on December 24\textsuperscript{th}, 1967, sparking not a few controversies (see \cite{G1}). The greatest writer of the Italian Literature prefers to write his treatises as dialogues, discourses. And to discourse is to run. Galileo is also the man who introduces the experimental method and the scientist of \emph{objects} like the telescope or the spheres he is said to have dropped from the Tower of Pisa. This lesson from Calvino is full of sought objects (Charlemagne's ring), proposed objects (the horse in a novella by Boccaccio), prepared objects (the drawing of Chuang-Tzu) and even non-existent objects (the book by Borges). The object aids the mental speed that Calvino assures it is a value. Even in this lesson, he emphasizes that 

\begin{quote}
Each value or virtue I chose as the subject for my lectures does not exclude its opposite. Implicit in my tribute to lightness was my respect for weight, and so this apologia for quickness does not presume to deny the pleasures of lingering. Literature has worked out various techniques for slowing down the course of time. I have already mentioned repetition, and now I will say a word
about digression. \cite{C2}.
\end{quote}


In this phrase appears the crucial element for speed: time. Confusing speed and velocity, it  appears obvious that time is crucial being velocity a function of time, derivative of space with respect to time. But speed and velocity are not the same thing in this lesson. The desire is for the speed of mental calculation, not the haste of the result. Even though the objects move quickly in the mathematician's mind, they remain static when the mathematician is not looking at them. Time applied to mathematical objects generates variations in the continuous as well as successions in the discrete. Calvino uses the words ``continuity''  but also ``discontinuity''  of time in this lesson. He talks about the succession of events like the rhymes of a poem. Let's revisit the objects of the Calvinian lesson in this temporal and mathematical key.

The ring sought by Charlemagne is magical, in the legend it connects all the characters searching for it. The mathematical object most similar to this ring is a vertex of a graph. Calvino clearly writes that the scheme of a story is a graph, well knowing the combinatorial topics that his friends from Oulipo are passionate about.

\begin{quote}
The events, however long they last, become puncti-form, connected by rectilinear segments, in a zigzag pattern that
suggests incessant motion. \cite{C2}.
\end{quote}

This zigzag alludes to a random walk that requires a pace similar to that of the horse in Boccaccio's novella, from which Madonna Oretta wants to dismount because there is speed without rhythm. 

\begin{quote}
In an
age when other fantastically speedy, widespread media are triumphing, and running the risk of
flattening all communication onto a single, homogeneous surface, the function of literature is
communication between things that are different simply because they are different, not blunting
but even sharpening the differences between them, following the true bent of written language. \cite{C2}.
\end{quote}

 Science too has always had this function: to highlight differences. For example, the notion of curvature can distinguish among different geometric worlds. It is precisely in this lesson that Calvino talks about geodesics, about shorter paths, about straight lines in a labyrinth but also about routes that connect distant points in space and time. He must have been influenced by the ideas of the theory of relativity, a noun that Calvino precisely uses in this lesson. 
 
 Borges pretends that his book has already been written, constructs imaginary libraries that open to infinity. Here is the true speed, the swift thought with which we imagine infinity compared to the effort of counting. 

On the other hand, the design of a crab made in 10 years with perfection in a single gesture can correspond in Mathematics to the distance between a question about infinity and the corresponding elegant answer. For example, Galileo lets his characters Salvati and Simplicio discuss the correspondence between whole numbers and square numbers. In a single gesture, two hundred years later, Cantor responds to the paradox with a paradise, that of transfinite numbers. 
 
 The squares and roots are in Salvati's phrase ``just as the products are called squares, the factors, that is, those that are multiplied, are called sides or roots; the other [numbers], which do not arise from numbers multiplied in themselves, are not squares otherwise. Hence, if I say, all numbers, including squares and non-squares, are more than just squares alone, I say something true: isn't it?'' \cite{Galileo}. The square and the root are taken up in Calvino's description of Borges' writing ``With Borges, a literature raised to the square and at the same time a literature as the square root extraction of itself is born.''  \cite{C2}.
 
 Ideas quickly draw a gigantic graph between disciplines. The random walk in knowledge seems confused, but it is rapid, like this paragraph, like the third American lecture.

\section{Euclidean Exactitude}\label{sec:Exactitude}

In common language, precision and Mathematics are closely intertwined, often exemplified by the expression: ``mathematical''  as a synonym for ``exact''. In the Treccani's Italian dictionary, the definition of ``mathematical''  as an adjective is as follows: ``having the character of objective, clear, and consequent precision typical of mathematical calculation'' \cite{D}. This definition implicitly identifies Mathematics with ``mathematical calculation''.

Anyone who studies Mathematics even briefly experiences a division between Mathematics and calculation. For instance, in Italy during the 1990s, a popular booklet of a hundred pages by Giuliano Spirito titled ``Matematica senza numeri''  (``Mathematics without numbers'') \cite{S} was widely read, even by high school students who struggled to relate the core subjects of the book (sets, logic, relations) with the concept of Mathematics learned in primary school. The takeaway message was that Mathematics encompassed more than just calculations, numbers, and polygons. Those who study Mathematics at a university level cannot simply equate the discipline with the science of calculations.

In reading Italo Calvino's third lecture titled ``Exactitude'', we find not only this dichotomy but also a phrase that could serve as an alternative definition of Mathematics:

\begin{quote}
My search for exactitude was
branching out in two directions: on the one side, the reduction of secondary events to abstract
patterns according to which one can carry out operations and demonstrate theorems; and on the
other, the effort made by words to present the tangible aspect of things as precisely as possible.
\end{quote}

The point is that if Mathematics is defined as the language of abstract schemes with which operations can be performed and theorems proved, then it is necessary to start from the definition of the entities themselves: its language must define every object of the discipline without ambiguity.

This is not a common feature in the language we usually use to communicate, and indeed, nuances are a powerful tool that poets use in their compositions. This type of ambiguity is the first enemy that a mathematician, in the work, must defeat so that a mathematical entity is the same for all who inhabit a certain abstract universe. For example, in Euclidean geometry, for everyone, a square is a polygon with four equal sides and four equal angles. The omission of any detail would lead to the construction of another figure: omitting only that the polygon has four equal angles leads to the construction of a rhombus, omitting only that the lengths of the sides are the same leads to the construction of a rectangle, and so on. Every word is essential to define the square exactly. At that point, we will need an example of a square, which we will find in an image, use in a building, or in a constructed object. The Lightness of the abstract entity, combined with the Exactitude of its definition, avoids emptiness. There is no opposition between Lightness and Exactitude.

We can conclude that in 1963, Italo Calvino turned to the use of science in the novel, for escaping a fear he reveals in this chapter of the Lectures:

\begin{quote}
It sometimes seems to me that a pestilence has struck the human race in its most distinctive
faculty—that is, the use of words. It is a plague afflicting language, revealing itself as a loss of
cognition and immediacy, an automatism that tends to level out all expression into the most
generic, anonymous, and abstract formulas, to dilute meanings, to blunt the edge of
expressiveness, extinguishing the spark that shoots out from the collision of words and new
circumstances. \cite{C2}.
\end{quote}

\begin{figure}[h!]
	\begin{center}
	\includegraphics[width=\textwidth]{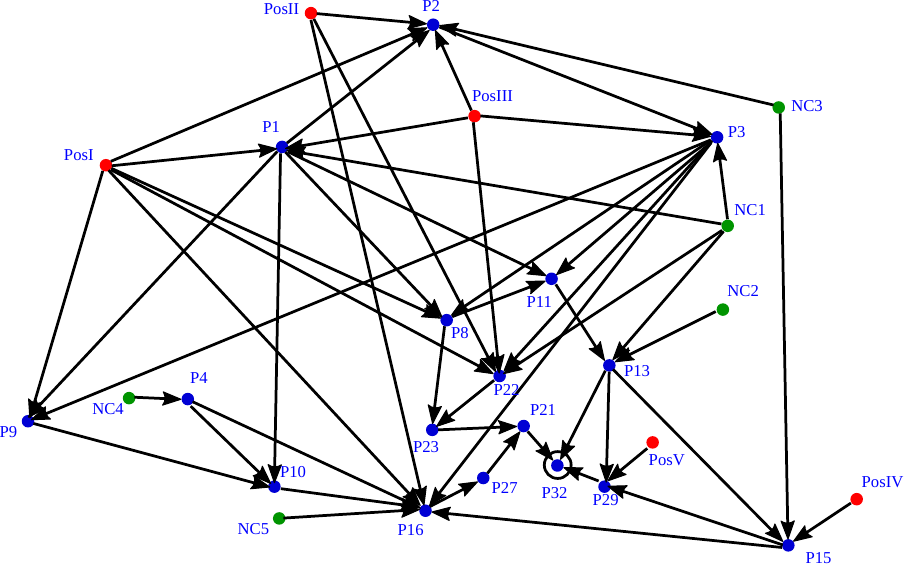}
\end{center}
\caption{All Propositions (P), Common Notions (CN) and Postulates (Pos) used in Euclid's elements to prove Proposition 32. An arrow going from a point A to a point B indicates that A is used to prove B.}\label{Fig1}
\end{figure}

Beyond abstract schemes, Calvino talks about theorem proofs. First, a theorem must be stated. The statement of a theorem requires the same precision used in definitions.
But a theorem has something more than a definition. A theorem is a statement in which, given some premises called hypotheses, some consequences called theses are asserted.

Consider Proposition 32 from Book I of Euclid's Elements: in his axiomatic system, the sum of the internal angles of a triangle equals a straight angle. A theorem is accompanied by its proof, which necessarily leads from the hypotheses to the theses. The difference between a philosophical discussion and the proof of a theorem is that in the theorem proof, every single step must be justified precisely. The order of the steps is crucial. In Figure \ref{Fig1} (taken from \cite{Lu1}), there is a directed graph with vertices representing all the propositions of Book I of the Elements up to the 32nd (numbered with P), the common notions (CN), and the five postulates (Pos). The edges of the graph indicate which result is necessary to deduce a certain proposition. We understand that the Fifth Postulate is necessary from Proposition 29 onwards, particularly to prove the above-mentioned Proposition 32. Moreover, the theorem on the sum of the internal angles of a triangle could not have been stated for example before Proposition 13, which is used in its proof. The deduction of a result from the previous ones is therefore similar to the process Calvino describes in this lecture on writing: the detail of the detail.

There are, however, several possible proofs for the same theorem: some more elegant and concise, others more cumbersome and convoluted.

Let's go back to to the non-contradictory relationship between Exactitude and Lightness. For Calvino, Exactitude is figuratively represented by the feather which in Egyptian tradition serves as a counterweight for the souls of the deceased reaching the afterlife. It's interesting the analogy of this point with the proof of a theorem which, in everyone's dreams, should be as beautiful and simple as a feather.

There's an interesting story about this. According to Paul Erd\"os, the most prolific mathematician of all time, God wrote all the most elegant mathematical proofs in a book, to which he referred as ``the Book''. To convey how important this pursuit of perfection, exactitude, lightness, and elegance was to him, let's quote what he said at a conference in 1985: ``you may not believe in God, but you must believe in the Book''. An attempt at ``the Book'' has been made by Aigner and Ziegler (see \cite{AZ}), but every mathematician must choose the proofs they would put in the book and ask themselves if their proof has the elegance to aspire some pages.

\section{Quantum Visibility}\label{sec:Visibility}

We are inhabitants of the third millennium, we regard the American Lectures also with the curiosity to understand if the values Calvino chose have really been preserved. Coming to \emph{Visibility}, we can only assert, meanwhile we receive an Instagram notification, that Visibility is the heart of our time, but here Calvino does not talk about shots, about selecting among them, he talks about how the images of the mind become visible.

\begin{quote}
I have in mind some possible pedagogy of the imagination that
would accustom us to control our own inner vision without suffocating it or letting it fall, on the
other hand, into confused, ephemeral daydreams, but would enable the images to crystallize into
a well-defined, memorable, and self-sufficient form, the \emph{icastic} form. \cite{C2}.
\end{quote}

Icastic means with effectiveness, with immediacy, with concreteness. An icastic description is not only to remember, it becomes the concept itself. Crystallizing an idea is a process that modern science had somehow rejected, but rereading Calvino's sentence we realize that every part of twentieth-century Physics has proposed itself through the crystallization of images. This is not a paradox, there are many crystals and put together they enhance their reverberation.

Icastic images are Newton's apple, Einstein's twins, Schr\"odinger's cat, Hilbert's beer mugs. These images were born with the theory, they are not just a didactic expedient for its explanation. Calvino wonders how images are born. We can ask ourselves how an entire theory is born. Instead of writing a treatise on neuroscience, to explain this \emph{memo}, he talks about fantasy:
\begin{itemize}
\item Fantasy is a place where it rains inside \cite{C2}
\item From where do the images rain into fantasy? \cite{C2}
\item Fantasy is a kind of electronic machine \cite{C2}
\end{itemize}
These three sentences serve to argue on whether imagination is or is not a tool of knowledge. In this work, we lean towards yes and make the example of Quantum Mechanics. This theory is never mentioned in the lessons, but Darconza and Bischi in \cite{BD} highlight the sentence found in the lesson on \emph{Multiplicity} ``Even before science had officially recognized the principle that observation somehow intervenes to modify the observed phenomenon'', Gadda knew that ``to know is to insert something into the real: it is, therefore, to deform''. In paragraph \ref{sec:Singularity} we will return to Gadda's influence on Calvino, here we are interested in highlighting that Calvino hints at the \emph{Heisenberg Uncertainty Principle}:
\begin{equation}
	\Delta x \Delta p \geq \frac{\hslash}{2}
\end{equation}

This principle has generated many images later used in science fiction stories: microscopic tunnels, collapsing particles, entanglements to be tattooed, up to teleportation. In these narratives, the conjugate variables are increasingly different: not only position-velocity, but for example rigor and clarity in popularization of science (see \cite{G2}) have correlated measures.

The uncertainty principle has entered the popular culture of this millennium precisely thanks to its Visibility despite initially appearing fanciful. Let's return to Calvino's lesson: from where have quantum ideas rained down? The very word ``quantum'' applied by Einstein to light is icastic. If there is a rain of images, then there must be, like for water, a cycle of ideas: they are the ideas of the scientists who preceded that image. Einstein's quantum would not have been there without Planck and Maxwell. Speaking of Maxwell, the formalization of the concept of field is also icastic. We don't need to explain the similarity between vectors and grass, it's immediate.

The icastic image of Schr\"odinger's cat to describe quantum entanglement is a mental experiment, so if fantasy rains inside literature to generate stories, it also rains in science to generate questions before explanations. For example, the cat led us to reflect on the macroscopic - microscopic transition. Schr\"odinger's cat is not a metaphor, but an analogy, and if Calvino in this lesson talks precisely about ``fields of analogies'', Lolli in \cite{L2} reminds us of the Langlands program to first connect Galois Groups, Algebraic Number Theory and then Harmonic Analysis up to Logic and Quantum Physics. Here again, the key concept is that of field, but we are interested in the attitude of those who are still today looking for a sort of unified theory even in Mathematics. The scientist is not so much the magician who activates the rain, but the scholar who is able to inhabit places where it rains and collect the rain: he inhabits laboratories without a roof of imposed choices, he opens windows on doubts and begins to observe the effect of the wind on the drops, he turns umbrellas upside down to collect the ideas of his predecessors. In a non-metaphorical sense, the imaginative scientist is attentive.

If for Calvino the \emph{imagination is the place of potential and hypothetical}, for a theoretical physicist hypotheses invent potentially real places. The effort in building a quantum computer today is not only linked to speed, efficiency, new applications, but because it must be a new electronic machine, overturning what Calvino said, there must be a kind of fantasy that actuates it.

\section{Algebraic Multiplicity and More}\label{sec:Multiplicity}

Multiplicity, in Mathematics, is an omnipresent concept, and is the second question one tries to answer when faced with a problem. The first question is the problem of existence which at first glance seems to be a naive question but is often anything but trivial. For example, an expert of differential equations can not avoid to ask whether the equation he uses to model reality has a solution or not. Indeed, a negative answer to this question implies that the model is completely useless from the point of view of modelization. But even a simple second-degree equation in the field of reals may or may not have solutions. But, assuming it does, how many are there? Practicing with decompositions and Delta, every high school student distinguishes algebraic equations with one or more solutions. Many also know different types of multiplicity (algebraic, geometric). The fundamental theorem of algebra due to Gauss in the complex field gives the number of these solutions.

For models that translate into differential equations, the same questions of existence and multiplicity depend on the meaning given to the word ``solution''. Therefore, the answer to the multiplicity of solutions depends on the spaces in which solutions are sought. To illustrate how these problems are of interest, it is useful to remember that one of the Millennium Prize Problems concerns the existence of solutions for equations describing fluid dynamics: the Navier-Stokes equations. In order to stay closer to the Calvinian text, we prefer to return to the example of algebraic equations starting from the following.

\begin{quote}
To show that in our own times literature is attempting to
realize this ancient desire to represent the multiplicity of relationships, both in effect and in
potentiality. \cite{C2}.
\end{quote}

Calvino's lesson on \emph{Multiplicity}, therefore, focuses on \emph{relationships}; we want to identify the world of relationships between various parts of Mathematics that determine the way in which mathematical knowledge progresses. In fact, Calvinian memos have a non-empty intersection: in the paragraph of \emph{Exactitude} Calvino reports the following quotation from ``The Man Without Qualities'' by Musil.

\begin{quote}
He mentions the fact that mathematical problems do not admit of a general solution, but that particular solutions, taken all together, can lead to a general solution. \cite{C2}.
\end{quote}

In Mathematics, there are very general results that are true milestones used daily to discover new things. In the romantic image of the mathematician (of the scientist in general), results come from brilliant and instantaneous intuitions; in the reconstruction that can be made of the history of Mathematics, this is not how discoveries are made. A result dealing with general problems is typically produced after years and years of work carried out independently by many people who, by tackling specific problems with very different techniques, go on to build bricks of knowledge. With many bricks available, it is possible to understand the possible relationships between these elements, that is, to arrange all the bricks coherently. Thus, the formulation of a general result is born. 

An example of this process is given by the history of solutions to third-degree equations. We are between the 15\textsuperscript{th} and 16\textsuperscript{th} centuries, and at that time, it was usual for mathematicians to compete in duels by solving mathematical problems. The consequences of a defeat could be dramatic, such as losing a job or going bankrupt; the consequences of a victory could be fame and glory until the next duel. It was therefore particularly important to find methods to solve a large number of problems in a general sense and to play this method time after time without revealing it. Initially, only third-degree equations of a very specific form were known to be solvable, and the first more general result is due to Scipione del Ferro, who managed to find a solution formula for third-degree equations missing the quadratic term. Later, another mathematician, Niccolò Fontana (nicknamed ``Tartaglia''  for his difficulty in speaking), discovered a way to solve all third-degree equations missing the linear term. It was then Gerolamo Cardano who discovered how to find the general solution of third-degree equations, showing that a cubic equation can always be reduced to an equation missing the quadratic term and thus, quoting Musil, he succeeded by combining individual solutions in a way to obtain the general solution. 

In modern Mathematics, the challenge is made by many more competitors on different themes. The lesson on \emph{Multiplicity} begins with a long quotation from Gadda and with the challenge posed by the encyclopedic novel. In the mathematical field, the encyclopedia of twentieth-century results lives in huge databases (ArXiv, ISI, Scopus, MathSciNet, Zentralblatt). The dissemination of results is itself a challenge. Carlo Emilio Gadda's quotation is used by Calvino to talk about a ``system of systems'' in literature. The mathematical counterpart is the intertwining of various disciplinary interests and the keywords of different scientific articles. According to Calvino, a second aspect of Multiplicity is related to what, according to Calvino, unites writers like Musil and Gadda, namely:

\begin{quote}
the inability to conclude. \cite{C2}.
\end{quote}

Imagine a situation where a group of mathematicians works for a certain amount of time on a problem until they find new, perhaps even non-trivial aspects. At that moment, the discussion begins on whether it is appropriate to spend part of the time to properly write down the result obtained and publish it or whether to focus all efforts on finding a more general result. On one hand, the advantage of having only sufficiently general results would allow a streamlining of the mathematical literature. On the other hand, the publication of minor results allows for easier understanding by the public and also a greater dissemination of what Musil calls individual solutions.

 This discussion hides, not too subtly, the fact that the definition of the conclusion of a scientific work is not unanimously shared. Learning to conclude a creative work is a complicated act so much so that ``Starting and finishing'' was among Calvino's memos (see Paragraph \ref{sec:BeginEnd}).
 
 Is it possible to say that a result is concluded? Probably this a Multiplicity of scientific relationships in
which the newer one does not prevail over the older one but, using a Newtonian image, allows us to see further by ``standing on the shoulders of giants''.

\section{The Missing Lesson: Consistency}\label{sec:Consistency}

No one can say what meaning Calvino would have given to the lesson on \emph{Consistency}; conjectures are made in \cite{P}. In \cite{L1} there is the autograph index with the word ``consistency'' written by Calvino with pencil compared to the previous five memos written with ink. Besides the syntactic aspect (Calvino would almost certainly have gone over the word with ink after writing the lecture), there is a semantic difference between consistency and the other memos: this strongest demand for a novel or a theorem had in the early twentieth century suffered a double attack. Well known to everyone is the literary Decadentism that has its scientific counterpart in the studies of Kurt G\"odel.

Consistency is a key word in logic: the property of a system of axioms for which it is not possible to deduce a contradiction from it.

Would a consistent literature be valuable in the sense that it gives us unambiguous outcomes of truth or falsehood? Mathematical consistency, on the other hand, refers to axioms and therefore, anyone wanting to answer such a question should consider the problem of ``axiomatic'' literature. Tracing back to literary, the axiomatic method would have been daring much more than doing combinatorial literature that Calvino had learned in Oulipo. In combinatorial literature, all the ingredients are already there and need to be mixed. In the axiomatic context, instead, one cannot predict a theorem just by looking at the axioms. \emph{Consistency} is linked not only to the idea of truth but also to that of \emph{proof}. Mathematicians use the terms \emph{consistency} and \emph{logical coherence} interchangeably. But in literary analyses, coherence is more related to the \emph{connection of the parts} of a whole, while consistency is related to the idea of \emph{solidity}. A mathematical proof is always solid and relative to what has been learned up to that moment, therefore what is consistent is always coherent. (see again Fig.\ \ref{Fig1} as example.)

To reinforce the hypothesis that Calvino would have consulted the logical-mathematical results, we recall the quote from Hofstadter \cite{H} present in the lecture on Visibility. The book \emph{G\"odel, Escher, Bach} by Douglas R.~Hofstadter had been published in Italy only the year before the writing of the Lectures, but it had struck Calvino, who was investigating imaginative processes. From this quote, we deduce that Calvino was familiar with the incompleteness theorems.

The word consistent appears almost as a constraint in Gödel's First Incompleteness Theorem. Let's rewrite it in a simplified manner: 

\emph{ ``Every consistent formalization of Mathematics that is sufficiently powerful to be able to axiomatize the elementary theory of natural numbers is capable of constructing a syntactically correct proposition that cannot be proved or refuted within the same system.''}

Incompleteness, therefore, is the possibility that paradoxes arise. Let us recall that Calvino cites undecidability in the lesson on Visibility. The paradoxical is of great interest to the literate. The G\"odel Theorem itself is based on the arithmetic transposition of the paradox ``this statement is not true''.

In the Preface to the Garzanti Italian edition of the American Lectures, Esther Calvino informs us that in the lecture on Consistency the author would have consulted Herman Melville, in particular the work \textit{Bartleby the Scrivener: A Story of Wall Street} \cite{Mel}. The character of this story is itself paradoxical: he refuses to perform tasks using the phrase ``I would prefer not to''. Denial is certainly not a value for the third millennium.

Mathematics overcame the impasse generated by G\"odel's Theorem similarly to how literature has overcome negatives for
not stop working (which happens in \emph{Bartleby}).
The second incompleteness theorem, also by G\"odel, asserts that
\emph{ ``Considering a mathematical theory sufficiently expressive to contain arithmetic; if it is
consistent, then it is not possible to prove the consistency of the theory within the
theory.''}

This too would have been a value for Calvino: distinguishing literature from the argument of literature. Perhaps this is why the novel is today extremely abundant as are the mathematical results: in the majority of cases the novel doesn't talk about itself, mathematicians don't do metatheory.
But there is also an incompleteness of the narrated compared to reality. Let's reread the incipit of Melville's story in which the first-person narrator warns the reader: ``it doesn't exist
material - I am convinced - to compose a complete and satisfactory biography of
this man.'' We, mathematicians, do not preface every work with the theorems of
incompleteness, but if we have to bet on why consistency is a value
for us, it is precisely by virtue of G\"odel's theorems: the mathematical creativity is contained
from the broad horizon of metatheory and the distant glimmer of paradoxes. What's beyond this mathematical and literary horizon? Calvino writes in the lesson
on \emph{Exactitude}:

\begin{quote}
But maybe this lack of substance is not to be found in images or in language alone, but in the
world itself. \cite{C2}.
\end{quote}

\section{The Incomplete Appendix: On the Beginning and the Ending}\label{sec:BeginEnd}

The lesson\emph{ ``On the Beginning and the Ending''} which appears in the appendix of \cite{C1} was probably meant to be an introduction to the six Harvard lectures. If these represent not only the values of literature but those of all knowledge, then these notes read as a response to the following questions: Why does one decide to write a certain story? When is a physical theory developed? How is the idea of a theorem generated? Calvino speaks of the world as a ``sum of information'' and language as something to be extracted from the languages developed by various disciplines. But immediately after, he writes that the beginnings lead us into a world different from the one that generated the idea. The mathematical parallel in this case could be the different geometries. For example, any version of the fifth axiom posed at the beginning of Euclid's masterpiece leads us into the new worlds of elliptic or hyperbolic geometry. Gauss's Theorema Egregium marks the beginning of differential geometry. The problem of the bridges of K\"onigsberg initiates the topological universes. The ancient gasket of Apollonius continues in fractal geometry. These are examples of preambles in Mathematics. Obviously, they are in relation, one with the others: differential geometry includes Euclidean and non-Euclidean; without the concept of topological dimension, there would be no description of fractal dimensions. As Calvino intimates with his literary examples, the preamble is a connection in the continuous connective tissue of human events. Calvino explains the choice of a single story (for us mathematicians, a problem to be addressed) with the words ``from the universal to the particular'', a story is ``isolated'' from others to be told. In the same way, in a scientific theory, an experiment or a problem becomes an object of great interest. Having spoken of topology, we cannot but give the example of the Poincaré Conjecture, which has led to great developments in the discipline of which Poincaré himself is considered the founder. Certainly, it is easier to discuss topology starting from Euler's result on graphs than from the concept of dimension in Poincaré. Here, not only those who study, but also those who disseminate Mathematics face the dilemma that Calvino expresses well in this appendix/introduction:

\begin{quote}
The moment of choice: we are given the possibility to say everything, in all possible ways; and we must come to say one thing, in a particular way. \cite{C2}.
\end{quote}

As in \cite{Lu2}, we draw a parallel between this quote and the axiom of choice stated in 1904 by Zermelo: ``Given a non-empty family of non-empty sets, there exists a function that assigns an element from each set of the family to it''. We are in set theory, we have an infinity of sets, and we want to call a representative from each set, a kind of symbolic reunion, is it possible? Do the other axioms of set theory already assure us of this? The history of this axiom, its acceptance, its consequences, and equivalences would make a beautiful story, but having to ``come to say one thing, in a particular way'' we limit ourselves to the concept of function. This has in its definition the notion of uniqueness (of the chosen representative) which in Calvino's manuscript appears as a challenge for each of us: ``in the multiplicity of possible experiences [...] the uniqueness of the days we have to live'' \cite{C1}. We might add the uniqueness of the theorems we have to prove and the stories we have to tell. The axiom does not provide an algorithm on how to choose; it is the creative zone of the mathematician and the writer. Calvino includes oblivion in the set of things to narrate (he writes precisely so, \emph{the set}), for mathematicians, however, the void is not a value of functions: one starts from the void to build for example, the natural numbers, but does not find the void at the conclusion. This asymmetry is also present in the literature. In this lesson, Calvino asks: ``Can we make symmetrical considerations for the ending as we have done for the beginning?'' It seems so because if the beginning is connected to human events, the ending (even of this incomplete lesson) is ``to continue telling'', i.e., the story returns to human events. However, everything is transformed, how to recognize the action of the story on these events?

 Returning to Mathematics, let's see the beginning and conclusion of an IFS fractal, i.e., obtained by iteration. If Calvino remembers that ``Dante ends the three parts of his poem with the word stars'', we can consider precisely in the arrangement of the stars a geometry of voids. What is the initial rule? What is the first verse of this astronomical poem? If Calvino hints at ``an indeterminate ending'', we can take that square (defined with Exactitude in Paragraph \ref{sec:Exactitude}) and start dividing it into nine equal squares, coloring the central square. Repeating this operation in the remaining eight squares, we get 9 colored squares (1 large, 8 small) and 64 uncolored small squares. We can start over, and the ``indeterminate ending'' after infinite steps will be Sierpinski's carpet square. It is one of the first fractals to appear on the mathematical scene after having made its presence known in finite steps in the floors of churches and castles. If Calvino writes of ``an ending that dissolves the illusion'', we could look at the Peano or Hilbert curve that seems like a thread but fills the square, therefore, it is two-dimensional, in contrast to the previous carpet that seems a rigid object of the plane but is actually one-dimensional. Finally, Calvino talks about an ending that ``questions the entire narration,'' and we can not avoid to think of the sensation one feels after the beautiful lessons on inductive topological dimension in which the topological dimension is defined in such a way that the previous topological dimension is assigned to the border of an object of the next topological dimension; the reassuring sensation contrasts with the appearance of the Mandelbrot set in the plane, which has a border of dimension 2. Clearly, the notion of dimension has changed, on the other hand, we have also written ``finally'' in the previous paragraph, we thus have what Calvino calls the ``cosmic ending'' or how every novel, just like every mathematical theory, is a universe in which the beginning are the axioms and the conclusion is just a passage to a new universe.

\section{Singularity: a Seventh Lesson?}\label{sec:Singularity}

A lesson on \emph{Singularity} would seem to be in apparent opposition to the lesson on \emph{Multiplicity}, somewhat like the relationship between \emph{Lightness} and \emph{Exactitude}; but this does not surprise us: as we have already seen, Calvino writes ``every value I choose does not intend to exclude the opposite value.'' \cite{C2}. From a literary perspective, this opposition is clearly described in \cite{N} also in relation to Gadda's writings: hence, to write, to begin a novel moves from the Multiplicity of stories that could be told to the single story that will become the subject. This choice will inevitably leave out an infinity of possibilities. But in hindsight, after the story is ready for the reader (or the Reader), the act of writing is a break from Multiplicity in favor of a Singularity, represented by the novel itself.

\begin{quote}
How is it possible to isolate a singular story if it implies other stories that cross and ``condition''  it, and these others still, extending to the entire universe? And if the universe cannot be contained in a story, how can one detach from this impossible story other stories that make complete sense? \cite{C2}.
\end{quote}

This way of seeing is not too far from how Singularity is discussed in the sciences. In both Mathematics and the natural sciences, at the point of Singularity, descriptive ambiguities are created, and thus, a large number of different coherent solutions can be proposed. Typically, most of them are not used in the surrounding context.

In Mathematics, the word singularity indicates a point in a certain space where a certain mathematical object is no longer well-defined as in all other points, making it special and interesting. In the more general scientific language, the word singularity indicates a typically extreme situation in which a fairly general theory stops working. Even in Mathematics and the natural sciences, Singularity and Multiplicity -- now understood in a different sense from Calvino's-- are two concepts inextricably linked. 

A typical example, which nowadays is perhaps the most well-known, is the theory of black holes. A black hole is an object where all the mass is concentrated at a point and has such a strong gravity that, beyond a limit called the ``event horizon'', not even light can escape. The equations of General Relativity that describe this phenomenon, the Einstein equations, predict a mathematical singularity at the point where all the mass is concentrated, and to date, there is no theory on how to describe what happens in such confined spaces and under such strong gravitational fields. This situation triggers an infinite multiplicity of possible descriptions among which, to date, none is experimentally verifiable. 

There are also \emph{singular phenomena} closer to our daily lives that are just less famous but no less interesting than black holes. An example is the phase transitions of substances and the existence of so-called critical points, where the usual phases of matter lose their meaning. At these special points, under very precise conditions of pressure, volume, and temperature, wonderful properties occur. In short, near a critical point, there are quantities that exhibit universal behavior, meaning they do not depend on all the details of the interaction between atoms and molecules, but only a tiny part of this information is retained, and therefore, very different substances exhibit the exact same behavior. Conceptually, here it is quite clear how the problem of singularity emerges: when dealing with the critical point, it must happen that the equations describing different substances become identical, and therefore, while it is possible to move towards the critical point, it is not possible to uniquely go back from the critical point. This is, in simple terms, one of the most beautiful ideas in the Theoretical Physics of the 20\textsuperscript{th}-century: the \emph{renormalization group}. 

An example of universal behavior is seen considering the singular limit of different equations such as the zero-viscosity limit of the hydrodynamics equations, the zero-range limit for the Schr\"odinger equation with potential, or the limits leading to the justifications of effective equations for the description of many-body systems.

Today, the study of these singular models is central both in Physics and Mathematics and represents not so much a value to carry into the future, but a starting point for the Mathematics of tomorrow.

\vspace{1cm}

\noindent
\textbf{Acknowledgements.}
This work for S.L. is part of the PNRR MUR project PE0000023-NQSTI NQSTI Spoke 9. M.G. was supported by the ERC StG project MaMBoQ, no. 80290 and the GNFM (National Group for Mathematical Physics). Both authors thank the INdAM for funding the outreach activity that led to this publication.

\bibliographystyle{plain}
\bibliography{biblio_chapter}

\begin{thebibliography}{00}
\bibitem{AZ} Aiger M., Ziegler G., \emph{Proofs from the Book}. Springer-Verlag, 2009.

\bibitem{BD} Bischi G.I., Darconza G., \emph{Calvino e la limpidezza della complessità. Tra Palomar e
Parisi.} Aras Edizioni 2023.

\bibitem{C} Cabibbo N., \emph{Tra fossili e particelle, la grande sfida della luce e dell’ombra. Tutto Scienze e Tecnologie}. La Stampa 27-01-2010. 

\bibitem{C1} Calvino I., \emph{Lezioni Americane, Sei proposte per il prossimo millennio.} Garzanti
1988 (Mondadori 1993).

\bibitem{C2} Calvino I., \emph{Six Memos for the Next Millennium}. Harvard University Press, 1988.

\bibitem{C3} Calvino I., \emph{Se una notte d’inverno un viaggiatore}. Einaudi 1979 (Mondadori 1994).

\bibitem{C4} Calvino I., \emph{Palomar}. Einaudi 1983 (Mondadori 1994).

\bibitem{D} Dizionario Treccani, https://www.treccani.it/vocabolario/matematico/ accessed 02.20.2024.

\bibitem{Galileo} Galilei G., \emph{Discorsi e dimostrazioni matematiche intorno a due nuove scienze}. Boringhieri (1958).

\bibitem{G1} Greco P., \emph{Calvino, dalla Terra alla Luna}. Il Bo Live, 21 January 2019.

\bibitem{G2} Greco P., \emph{Quale comunicazione della scienza per i paesi emergenti}. J. Science
Comm. 4(3), 1-6 (2005).

\bibitem{H} Hofstadter D. R., \emph{G\"odel, Escher, Bach: un’eterna Ghirlanda Brillante}. Adelphi 1984.

	\bibitem{L1} Lolli G., \emph{Discorso sulla matematica. Una rilettura delle Lezioni americane di Italo Calvino}. Bollati Boringhieri 2011.
	
\bibitem{L2} Lolli G., \emph{Molteplicità potenziale e creatività al tempo del computer: un matematico
del 2000 legge Calvino}. Californian Italian Studies, Volume 12, 2023.

\bibitem{Lu1} Lucente S., \emph{Elementi dei grafi e grafi negli Elementi}. Paginaria Edizioni 2020.

\bibitem{Lu2} Lucente S., \emph{Il momento della scelta}. Prisma n. 55 (2023) pp. 40-41.

\bibitem{Mel} Melville H., \emph{Bartleby, the Scrivener. A story of Wall-Street}. Melville House Publishing 2004.

\bibitem{N} Niccolai S., \emph{Il “Sistema del mondo”. Calvino e l’eredità di Gadda}. Italianistica: Rivista
di letteratura Italiana, Vol. 34 N. 3 pp 29-43.

\bibitem{P} Piacentini A., \emph{Consistenza, l’inesplorata sesta lezione di Calvino}. Edizioni Progetto
Cultura 2004.

\bibitem{S} Spirito G., \emph{Matematica senza numeri}. Newton Compton Editori 1995.
\end{thebibliography}

\end{document}